\documentclass[11pt]{amsart}

\usepackage{commandearticle_new}
\usepackage{slashbox}
\usepackage{pgf,tikz}
\usetikzlibrary{arrows}

\makeatletter
\def\subsection{\@startsection{subsection}{2}%
	\z@{.5\linespacing\@plus.7\linespacing}{.25\linespacing}%
	{\normalfont\bfseries}}
\def\subsubsection{\@startsection{subsubsection}{3}%
	\z@{.5\linespacing\@plus.7\linespacing}{.25\linespacing}%
	{\normalfont\itshape}}
\makeatother
\graphicspath{{Images/}}
\newcommand{\grandtraittop}{\rule{\textwidth}{0.25em}\par\vskip0.5\baselineskip}
\newcommand{\grandtraitbottom}{\par\vskip0.05\baselineskip\rule{\textwidth}{0.25em}}

\everymath{\displaystyle}
\begin{document}

\begin{abstract}
This paper studies the stability properties of a two dimensional relative velocity scheme for the Navier-Stokes equations. This scheme inspired by the cascaded scheme has the particularity to relax in a frame moving with a velocity field function of space and time. Its stability is studied first in a linear context then on the non linear test case of the Kelvin-Helmholtz instability. The link with the choice of the moments is put in evidence. The set of moments of the cascaded scheme improves the stability of the d'Humières scheme for small viscosities. On the contrary, a relative velocity scheme with the usual set of moments deteriorates the stability.
\end{abstract}

\title[On the stability of a relative velocity lattice Boltzmann scheme]{On the stability of a relative velocity\\ lattice Boltzmann scheme \\for compressible Navier-Stokes equations}

\author{François Dubois}
\address[François Dubois]{CNAM Paris, d\'epartement de math\'ematiques, Univ Paris Sud, Laboratoire de math\'ematiques, UMR 8628, Orsay, F-91405} 
\email{Francois.Dubois@math.u-psud.fr}

\author{Tony Février}
\address[Tony Février]{Univ Paris Sud, Laboratoire de math\'ematiques, UMR 8628, Orsay, F-91405, Orsay, F-91405}
\email{Tony.Février@math.u-psud.fr}

\author{Benjamin Graille}
\address[Benjamin Graille]{Univ Paris Sud, Laboratoire de math\'ematiques, UMR 8628, Orsay, F-91405, CNRS, Orsay, F-91405}
\email{Benjamin.Graille@math.u-psud.fr}

\date{\today}

\maketitle

\section*{Introduction}

The lattice Boltzmann schemes have been successfully used for the simulation of the compressible Navier-Stokes equations in two or three dimensions \cite{Benzi:1992:0,Chen:1992:0,Qian:1992:0,dHu:1992:0}. This method aims to mimic the microscopic behaviour in order to simulate some macroscopic problems. The algorithm consists in evaluating some particle distributions. The particles, moving from node to node of a lattice, undergo a phase of collision and a phase of transport. Different collision operators have been proposed for the simulation of the Navier-Stokes equations. The simplest one is the single relaxation time operator \cite{BGK:1954:0,Benzi:1992:0,Chen:1992:0,Qian:1992:0} also called BGK. An alternative called the multiple relaxation times (MRT) operator \cite{dHu:1992:0,LalLuo:2000:0} has been proposed. During the collision, some moments, linear combinations of the particle distributions, relax towards the equilibrium with a priori different velocities. It contains the particularity to offer more degrees of freedom to fix the different parameters as the viscosities. The multiple relaxation times approach is thus more flexible than the BGK. Both schemes have been well studied particularly in terms of stability \cite{LalLuo:2000:0}. They still encounter some instability features as the viscosities tend to zero that limits high Reynolds number simulations.

In 2006, a cascaded scheme improving the stability for low viscosities has been presented \cite{Geier:2006:0}. Its relaxation occurs in a frame moving with the fluid velocity. To understand the positive features of this scheme, a general notion of relative velocity schemes was defined \cite{Fev:2014:0}. Their relaxation is made for a set of moments depending on a velocity field function of space and time that is the velocity fluid for the cascaded scheme \cite{Fev:2014:0} and zero for the d'Humières scheme \cite{dHu:1992:0}. These relative velocity schemes are not restricted to the simulation of the Navier-Stokes equations: they are defined for an arbitrary number of conservation laws.
Their consistency has already been studied for one and two conservation laws \cite{Fev:2014:0,Fev:2014:1,Fev:2014:2} but the same reasoning holds for an arbitrary number of conservation laws. 

The purpose of this contribution is to present some numerical stability results of the two dimensional nine velocities ($\ddqn$) relative velocity scheme for the compressible Navier-Stokes equations. We want to characterize the influence of the relative velocity and the link with the moments choice: the polynomials defining the moments of the cascaded scheme are different from the usual ones and may act on the stability. In a first part, we recall the basis of the relative velocity schemes. We then present the relative velocity $\ddqn$ we are interested in. The second part exhibits the results of stability, first in a linear context ($L^2$ von Neumann notion) and then for a non linear test case, the Kelvin-Helmholtz instability. It puts in evidence the link between the relative velocity, the choice of the polynomials defining the moments and the stability. 

\label{}

\section{Description of the scheme}
\label{sec:S2}
We first introduce the relative velocity scheme for an arbitrary number of dimensions and velocities. We then particularize it to the case of two dimensions and nine velocities.
\subsection{The relative velocity $\ddqq$ scheme}
\label{sec:S21}

This section presents the derivation of the relative velocity lattice Boltzmann schemes introduced in \cite{Fev:2014:0} and inspired by the cascaded scheme \cite{Geier:2006:0}.
Let $\lattice$ be a cartesian lattice in $d$ dimensions with a typical mesh size $\dx$. The time step $\dt$ is linked to the space step by the acoustic scaling $\dt=\dx/\lambda$ for $\lambda\in\R$ the velocity scale. We introduce $\vectV=(\vj[0],\ldots,\vj[q-1])$ a set of $q$ velocities of $\R^d$. This defines the scheme called $\ddqq$. We assume that for each node $\vectx$ of the lattice $\lattice$, and  each $\vj$ in $\vectV$, the point $\vectx+\vj \dt$ is still a node of $\lattice$.
The $\ddqq$ scheme computes a particle distribution 
$\vectf = (\fj[0],\ldots,\fj[q-1])$
on the lattice $\lattice$ at discrete values of time. An iteration of the scheme consists in two phases: the relaxation that is non linear and local in space, and the linear transport solved exactly by a characteristic method.\\

The relaxation phase reads more easily in a moments basis using the d'Humières framework \cite{dHu:1992:0}. A velocity field $\utilde(\vectx,t)$ that depends on space and time being given, we define the matrix of moments $\MatMu$ by
\begin{equation*}\label{eq:MatMu}
 \Miju[kj] = \Pk(\vj-\utilde), \qquad 0\leq k,j\leq q{-}1,
\end{equation*}
where $(\Pk[0],\ldots,\Pk[q-1])$ are some polynomials of $\R[\vars{X}{1},\ldots,\vars{X}{d}]$. 
This matrix of moments, supposed to be invertible, defines the moments $\vectmu=(\mku[0],\ldots,\mku[q-1])$ by the relation
\begin{equation}\label{eq:ftomu}
\vectmu = \MatMu \; \vectf,
\end{equation}
where $\mku$ is the \kieme moment.\\
 
The relative velocity schemes use a diagonal relaxation phase in the shifted moments basis 
\begin{equation}\label{eq:relaxationu}
\mkue = \mku{+}\sk (\mkueq{-}\mku), \qquad 0\leq k\leq q{-}1,
\end{equation}
where $\mkueq$ is the \kieme moment at equilibrium and $\sk$, the relaxation parameter associated with the \kieme moment for $0\leq k\leq q{-}1$. Some of these moments are conserved by the relaxation: they are associated with relaxation parameters equal to zero.
The equilibrium derives from a vector of distribution functions at the equilibrium $\vectfeq$ independent of the velocity field $\utilde$,  only dependent of the conserved moments.
\begin{equation}\label{eq:mueq}
 \vectmequ = \MatMu \vectfeq.
\end{equation}
The inverse of the matrix of moments is used to return to the distributions
\begin{equation}\label{eq:mtof}
\vectfe = \MatMinvu\vectmue.
\end{equation}
The transport phase spreads the particle distributions on the neighbouring nodes
\begin{equation*}\label{eq:transport}
 \fj(\vectx,t+\dt) = \fje(\vectx-\vj\dt,t), \qquad 0\leq j\leq q{-}1.
\end{equation*}
This framework embeddes the d'Humières scheme for $\utilde$ equal to $\vectz$ and the cascaded scheme for $\utilde$ equal to the fluid velocity and a particular set of moments \cite{Fev:2014:0}. In the following, when $\utilde$ is specified, we call the associated relative velocity scheme the scheme relative to $\utilde$.

\subsection{The study framework: the relative velocity $\ddqn$ scheme}\label{sub:ddqn}

The purpose of this section is to introduce the scheme whose stability properties are investigated: the relative velocity $\ddqn$ scheme with two conservation laws on the density and the momentum
\begin{equation*}\label{eq::rhoq}
\rho=\sum_j \fj,\quad q^{\alpha}=\sum_j v_j^{\alpha}\fj,\quad 1\leq\alpha\leq d.
\end{equation*} for the compressible Navier-Stokes equations. We expose its features and the different degrees of freedom used to check its stability. We put a particular attention on the definition of the moments.\\

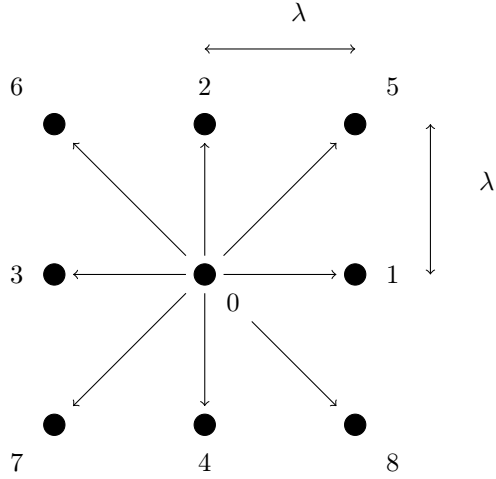
\begin{figure}
\begin{center}
\begin{tikzpicture}[scale=2.5,inner sep = 1.mm] 
		\tikzstyle{every node}=[font=\small]
		\draw[<->](0,1.2) -- (0.8,1.2);\draw[<->](1.2,0) -- (1.2,0.8);
		\draw[<-](-0.7,0) -- (-0.1,0); \draw[->] (0.1,0) -- (0.7,0.);
		\draw[<-] (0,-0.7) -- (0,-0.1); \draw[->] (0,0.1) -- (0,0.7);
		\draw[<-] (-0.7,-0.7) -- (-0.1,-0.1); \draw[->] (0.1,0.1) -- (0.7,0.7);
		\draw[<-] (-0.7,0.7) -- (-0.1,0.1); \draw[->] (0.25,-0.25) -- (0.7,-0.7);
	         \node[circle,draw,fill=black] at (0,0) {};\node at (-1,0) {3};\node at (1,0) {1};\node at (0,1.) {2};
	         \node at (0.15,-0.15) {0};\node at (0,-1) {4};\node at (1,1) {5};\node at (-1,1) {6};\node at (-1,-1) {7};\node at (1,-1) {8}; \node at (0.5,1.4) {$\lambda$};\node at (1.5,0.5) {$\lambda$};
	         \node[circle,draw,fill=black] at (0,-0.8) {};\node[circle,draw,fill=black] at (0,0.8) {};\node[circle,draw,fill=black] at (-0.8,0) {};\node[circle,draw,fill=black] at (0.8,0) {};\node[circle,draw,fill=black] at (0.8,-0.8) {};\node[circle,draw,fill=black] at (-0.8,-0.8) {};\node[circle,draw,fill=black] at (-0.8,0.8) {};\node[circle,draw,fill=black] at (0.8,0.8) {};
\end{tikzpicture}
\end{center}
\caption{The $\ddqn$ velocities.}
\label{fig:ddqn}
\end{figure}

For this two-dimensional scheme, nine velocities are involved: they are defined by
\begin{equation*}
\vectv=\{(0,0),(\lambda,0),(0,\lambda),(-\lambda,0),(0,-\lambda),(\lambda,\lambda),(-\lambda,\lambda),(-\lambda,-\lambda),(\lambda,-\lambda)\},
\end{equation*}
 with $\lambda\in\R$ the velocity scale. These velocities are also represented on the figure \ref{fig:ddqn}.\\
  
We need to deal with the set of the moments and the equilibrium to completely characterize the scheme. Historically the $\ddqn$ scheme has been mainly used with the following set of moments
\begin{equation}\label{eq:polddqn}
1,X,Y,X^2+Y^2,X^2-Y^2,XY,X(X^2+Y^2),Y(X^2+Y^2),(X^2+Y^2)^2,
\end{equation}
or its orthogonalized analogue for the simulation of the Navier-Stokes equations \cite{LalLuo:2000:0}. They have been chosen because of their physical meaning: they involve the density, the momentum, the energy, the diagonal and off-diagonal components of the stress tensor, the heat flux and the square of the energy. Nevertheless, the $\ddqn$ cascaded scheme \cite{Geier:2006:0}, that seems to improve the stability at low viscosities, has brought to light an other set of moments given by 
\begin{equation}\label{eq:polgeier}
1,X,Y,X^2+Y^2,X^2-Y^2,XY,XY^2,YX^2,X^2Y^2.
\end{equation}
This scheme has been written as a relative velocity scheme for the moments (\ref{eq:polgeier}) and the fluid velocity \cite{Fev:2014:0}. The relaxation of (\ref{eq:polddqn}) and (\ref{eq:polgeier}) are equivalent in the d'Humières framework corresponding to $\utilde=\vectz$ (section \ref{sub:stlindH}). However, this is not true any more when $\utilde$ is different from $\vectz$. This point naturally leads to the question: does the set of moments have an influence on the stability properties of the relative velocity scheme? Giving some experimental rudiments of an answer is the purpose of this study.\\

That's why we introduce two sets of moments tuned by a parameter $\alpha\in\R$:
\begin{equation}\label{eq:polal}
1,X,Y,X^2+Y^2,X^2-Y^2,XY,X(\alpha X^2+Y^2),Y(X^2+\alpha Y^2),\frac{\alpha}{2}(X^4+Y^4)+X^2Y^2,
\end{equation}
and
\begin{equation}\label{eq:polal2}
1,X,Y,X^2+Y^2,X^2-Y^2,XY,XY^2+\alpha(X^2+Y^2),YX^2+\alpha(X^2+Y^2),X^2Y^2.
\end{equation}

The moments (\ref{eq:polal}) generalize those of the $\ddqn$ cascaded scheme corresponding to $\alpha=0$ \cite{Fev:2014:0} given by (\ref{eq:polgeier}) and the ones associated with $\alpha=1$ defined by (\ref{eq:polddqn}). The introduction of $\alpha$ results from the will not to restrict the study to two sets of moments. This allows also to understand the impact on the stability of the $X^3$ component when $\alpha$ moves from $0$ to $1$.
The choice of the moments (\ref{eq:polal2}), even if it seems strange because it mixes some second and third order polynomials, improves the understanding of the differences of stability between (\ref{eq:polddqn}) and (\ref{eq:polgeier}) for the relative velocity  $\ddqn$ scheme. Taking $\alpha=0$ also recovers the cascaded moments.\\

The equilibrium may also have an influence on the stability. That's why, denoting $\vectu=\vectq/\rho$ the fluid velocity, we introduce
\begin{equation}\label{eq:eqqian}
\fjeq(\rho,\vectu)=\rho\vars{\omega}{j}\Big(1+\frac{\vectu.\vj}{\co[2]}+\frac{(\vectu.\vj)^2}{2\co[4]}-\frac{|\vectu|^2}{2\co[2]}\Big),\quad 0\leq j\leq8,
\end{equation}
and 
\begin{multline}\label{eq:eqgeier}
\fjeq(\rho,\vectu)=\rho\vars{\omega}{j}\Big(1+\frac{\vectu.\vj}{\co[2]}+\frac{(\vectu.\vj)^2}{2\co[4]}-\frac{|\vectu|^2}{2\co[2]}\\+\frac{(\vectu.\vj)^3}{6\co[6]}-\frac{|\vectu|^2(\vectu.\vj)}{2\co[4]}+\frac{d_j^{}(\ux)^2(\uy)^2}{\co[4]}\Big),\quad 0\leq j\leq8,
\end{multline}
where $d_0^{}=-1/4$, $d_j^{}=1/2$, $j=1,\ldots,4$, $d_j^{}=-1$, $j=5,\ldots,8$, 
respectively corresponding to the second order truncated equilibrium \cite{Qian:1992:0} and to the product equilibrium used for the $\ddqn$ cascaded scheme \cite{Geier:2006:1}. The product equilibrium corresponds to the fourth order truncation of the maxwellian equilibrium. Both equilibria allow to simulate the compressible Navier-Stokes equations whatever the velocity field $\utilde$. Indeed, the second order equivalent equations of the relative velocity schemes are independent of $\utilde$ \cite{Fev:2014:0} and the Navier-Stokes equations are recovered by the $\ddqn$ d'Humières scheme ($\utilde=\vectz$) at the second order for small Mach numbers \cite{Dub:2008:0}. Let's note that the simulations of the section \ref{se:KHstab} have been also made for the incompressible analogue of (\ref{eq:eqqian}) and (\ref{eq:eqgeier}) used in \cite{LalLuo:2000:0}: same trends as those presented in the section \ref{se:KHstab} are obtained.\\

We choose to work with several two relaxation times schemes (TRT) to understand the role of each polynomial of the moments: the one given by 
\begin{equation}\label{eq:trt1}
\vects=(\sk[e],\sk[\nu],\sk[\nu],\sk[e],\sk[e],\sk[e]),
\end{equation}
called TRT$_1$ and the one given by
 \begin{equation}\label{eq:trt2}
 \vects=(\sk[e],\sk[e],\sk[e],\sk[p],\sk[p],\sk[e]),
  \end{equation}
  called TRT$_2$ where $\sk[e],\sk[\nu],\sk[p]\in\R$. Note that the TRT$_1$ and the TRT$_2$ differ from the TRT schemes defined in \cite{Ginz:2010:0} and based on the symmetry of the lattice. If all the relaxation parameters are identical, we recover the BGK scheme \cite{dHu:1992:0}.\\
  
Four degrees of freedom are tunable in this section: the moments, the vector of the relaxation parameters $\vects$, the velocity field $\utilde$ and the equilibrium. The link between these parameters and their influence on the stability is studied in the following.

\section{Experimental study of linear stability}\label{se:num} 

In this section, we study the linear von Neumann $L^2$ stability of the relative velocity $\ddqn$ scheme defined in the section \ref{sub:ddqn}. The influence of the moments according to the choice of the velocity field $\utilde$ is the keypoint of the section. Our first interest goes to the moments (\ref{eq:polgeier}) and (\ref{eq:polddqn}), respectively corresponding to $\alpha=0$ and $\alpha=1$ in (\ref{eq:polal}), because they are usually chosen by the community \cite{dHu:1992:0,LalLuo:2000:0,Geier:2006:0,Wang:2013:0}. The first subsection compares those two sets according to the velocity field parameter. We show that taking $\utilde=\vectu$, the velocity of the fluid, improves the stability if the moments (\ref{eq:polgeier}) are chosen, deteriorates it if the set (\ref{eq:polddqn}) is taken.
The second subsection answers the following question: what is the better choice of moments (of $\alpha$) in terms of stability? A range of $\alpha$ is proposed and $\alpha=0$ is showed to be the most stable choice.

\subsection{Methodology: the von Neumann stability}

The study of the relative velocity $\ddqn$ scheme is based on the $L^2$ von Neumann stability. This notion being adapted to linear contexts, we linearize the equilibria (\ref{eq:eqqian}) and (\ref{eq:eqgeier}) around a velocity $\vectV=|\vectV|e^{i\theta}\in\R^2$, $\theta\in\R$. Thus there exists a matrix $\MatE$ so that  $$\vectfeq=\MatE\vectf.$$ Using (\ref{eq:ftomu},\ref{eq:relaxationu},\ref{eq:mueq},\ref{eq:mtof}) the linearized relaxation phase of the relative velocity schemes reads
\begin{equation*}\label{eq:rel}
\vectfe=(\MatI+\MatMu^{-1}\MatD\MatMu(\MatE-\MatI))\vectf,
\end{equation*}
where $\MatD={\rm diag}(\vects)$ is the diagonal matrix of the relaxation parameters. This expression holds for each node $\vectx$ of the lattice, the relaxation being local in space. One can deduce the expression of the distribution after an iteration thanks to the transport phase
\begin{equation*}\label{eq:iter}
\fj(\vectx,t+\dt)=[(\MatI+\MatMu^{-1}\MatD\MatMu(\MatE-\MatI))\vectf]_j^{}(\vectx-\vj\dt,t),\quad\vectx\in\lattice,\quad t\in\R.
\end{equation*}
In the Fourier space, the transport operator becomes local in space and is represented by the diagonal matrix $\MatA$ whose diagonal components are given by $e^{i\dt\vectk.\vj}$, $0\leq j\leq8$.
 We can then define the amplification matrix $\MatL(\utilde)=\MatL(\utilde,\vectV,\vectk,\vects,\alpha)=\MatA(\MatI+\MatMu^{-1}\MatD\MatMu(\MatE-\MatI))$, for $\vectk, \vectV, \utilde\in\R^2$, $\vects\in\R^9$, $\alpha\in\R$, $\vectV\in\R^2$ characterizing a time iteration of the scheme in the Fourier space
$$\widehat{\vectf}(\vectk,t+\dt)=\MatL(\utilde)\widehat{\vectf}(\vectk,t),\quad t\in\R,$$
where $\widehat{\vectf}$ is the Fourier transform of $\vectf$.
We want to determine the quantity
\begin{equation}\label{eq:maxVrho}
\max\{|\vectV|,~ \underset{\vectk\in\R^2}{\max}~r(\MatL(\utilde))\leq1\},
\end{equation}
for some parameters $\vects$, $\utilde$, $\alpha$, a direction of linearization $\theta\in\R$ and $r(\MatL(\utilde))$ the spectral radius of $\MatL(\utilde)$. It characterizes the set of the linearization velocities $\vectV$ for which the scheme verifies the necessary condition of $L^2$ stability $\underset{\vectk\in\R^2}{\max}~r(\MatL(\utilde))\leq1$.\\

\subsection{Comparison between the d'Humières scheme and the scheme relative to the linearization velocity}\label{sub:stlindH}

We show that the schemes relative to $\utilde=\vectV$ can improve or deteriorate the linear stability compared to the d'Humières scheme. The stability behaviour depends strongly on the choice of the moments.\\

We compare the schemes relative to $\utilde=\vectz$ and $\utilde=\vectV$ for the two sets of moments (\ref{eq:polddqn}) and (\ref{eq:polgeier}): we have $\alpha=0,1$ in (\ref{eq:polal}). We here restrict to the second order truncated equilibrium (\ref{eq:eqqian}) linearized around $\vectV$. The variable of comparison is the largest stable velocity $\vectV$ (\ref{eq:maxVrho}) for a linearization direction $\theta$ equal to $0$. We choose to deal with the TRT$_1$  (\ref{eq:trt1}) for $\sk[e]=2-2^{-m}$ and $\sk[\nu]=2-2^{-n}$ where $m,n,\in\N$, $0\leq m,n\leq7$. The parameters $\sk[e]$ and $\sk[\nu]$ respectively tune the bulk and the shear viscosities of the Navier-Stokes equations. This choice of parameters allows to study the zero viscosity limit by increasing $m$ or/and $n$. The table \ref{table:Mahum} deals with the d'Humières scheme for both sets of moments, the values for those two sets being identical.
The tables \ref{table:Mageier} and \ref{table:Mageier2} give analogous results for the scheme relative to $\utilde=\vectV$: they correspond respectively to the moments with $\alpha=0$ (\ref{eq:polgeier}) and $\alpha=1$ (\ref{eq:polddqn}).\\

\begin{table}
\centering\small
\grandtraittop
\begin{tabular}{@{}p{3.25cm}p{1cm}p{1cm}p{1cm}p{1.cm}p{1cm}p{1cm}p{1cm}p{1cm}@{}}
$n\qquad m$& 0&1&2&3&4&5&6&7\\ \hline
     0 & 0.42 & 0.41& 0.34&0.26&0.20&0.15&0.11&0.08 \\ 
     1 & 0.42 &0.41&0.36&0.30&0.23&0.18&0.13&0.09\\ 
     2 &0.31&0.34 & 0.34&0.32&0.28&0.23&0.17&0.13 \\ 
     3 &0.21&0.28&0.32&0.30&0.25&0.22&0.18&0.15  \\ 
     4 &0.14&0.21&0.28&0.26&0.22&0.18&0.16&0.13  \\ 
     5 &0.10&0.16&0.22&0.23&0.20&0.17&0.13&0.11  \\ 
     6 &0.07&0.12&0.17&0.20&0.18&0.16&0.12&0.11  \\ 
     7 &0.05&0.08&0.12&0.17&0.16&0.15&0.11&0.11  \\ 
\end{tabular}
\grandtraitbottom
\vspace{0.2cm}
 \caption{Highest stable $\vectV=(\Vx,0)$ in $\lambda$ units for the d'Humières scheme ($\utilde=\vectz$), with $\alpha=0$ or $\alpha=1$, of equilibrium (\ref{eq:eqqian}). $\sk[e]=2-2^{-m}$ and $\sk[\nu]=2-2^{-n}$.}
 \label{table:Mahum}
\end{table}

\begin{table}
\centering\small
\grandtraittop
\begin{tabular}{@{}p{3.25cm}p{1cm}p{1cm}p{1cm}p{1.cm}p{1cm}p{1cm}p{1cm}p{1cm}@{}}
$n\qquad m$& 0&1&2&3&4&5&6&7\\ \hline
    0 & 0.42 & 0.42& 0.40&0.37&0.33&0.28&0.24&0.21 \\
     1 & 0.42 &0.41&0.37&0.34&0.33&0.30&0.27&0.23\\ 
     2 &0.42&0.36 & 0.34&0.33&0.29&0.25&0.22&0.19 \\ 
     3 &0.36&0.34&0.33 &0.30&0.25&0.21&0.18&0.16  \\ 
     4 &0.32&0.32&0.30 &0.26&0.22&0.18&0.16&0.13  \\ 
     5 &0.29&0.30&0.27 &0.24&0.21&0.17&0.13&0.11  \\ 
     6 &0.26&0.28&0.24&0.21&0.18&0.16&0.12&0.11  \\
     7 &0.23&0.26&0.22&0.19&0.16&0.15&0.11&0.11  \\
\end{tabular}
\grandtraitbottom
\vspace{0.2cm}
 \caption{Highest stable $\vectV=(\Vx,0)$ in $\lambda$ units for the scheme relative to $\utilde=\vectV$ with $\alpha=0$ of equilibrium (\ref{eq:eqqian}). $\sk[e]=2-2^{-m}$ and $\sk[\nu]=2-2^{-n}$.}
 \label{table:Mageier}
\end{table}

\begin{table}
\centering\small
\grandtraittop
\begin{tabular}{@{}p{3.25cm}p{1cm}p{1cm}p{1cm}p{1.cm}p{1cm}p{1cm}p{1cm}p{1cm}@{}}
$n\qquad m$& 0&1&2&3&4&5&6&7\\ \hline
    0 & 0.42 & 0.42& 0.27&0.18&0.12&0.08&0.06&0.04 \\ 
     1 & 0.42 &0.41&0.37&0.31&0.20&0.13&0.09&0.06\\ 
     2 &0.24&0.36 & 0.34&0.32&0.27&0.21&0.14&0.09 \\ 
     3 &0.15&0.29&0.33 &0.30&0.25&0.21&0.18&0.14  \\ 
     4 &0.10&0.19&0.29  &0.26&0.22&0.18&0.16&0.13  \\ 
     5 &0.07&0.12&0.20 &0.29&0.21&0.17&0.13&0.11  \\ 
     6 &0.05&0.09&0.14&0.20&0.18&0.16&0.12&0.11  \\ 
     7 &0.03&0.06&0.09&0.14&0.16&0.15&0.11&0.11  \\ 
\end{tabular}
\grandtraitbottom
\vspace{0.2cm}
 \caption{Highest stable $\vectV=(\Vx,0)$ in $\lambda$ units  for the scheme relative to $\utilde=\vectV$ with $\alpha=1$ of equilibrium (\ref{eq:eqqian}). $\sk[e]=2-2^{-m}$ and $\sk[\nu]=2-2^{-n}$.}
 \label{table:Mageier2}
\end{table}

We notice the importance of the choice of the moments for the schemes relative to $\utilde=\vectV$: stability areas are the biggest for $\alpha=0$ (\ref{eq:polgeier}) (table \ref{table:Mageier})  and the smallest for $\alpha=1$ (\ref{eq:polddqn}) (table \ref{table:Mageier2})  whatever the choice of $\vects$. The d'Humières scheme ($\utilde=\vectz$, table \ref{table:Mahum}) has smaller stability areas than the scheme relative to $\utilde=\vectV$ with $\alpha=0$ and bigger than the one with $\alpha=1$.\\

The scheme relative to $\utilde=\vectV$ with $\alpha=0$ provides the most important gain compared to the d'Humières scheme when $\sk[e]$ or $\sk[\nu]$ is close to $2$ and the other is far from $2$ (for one small and one large viscosity). Instead these areas are the most deteriorated when $\alpha$ is equal to $1$. When $\sk[e]$ and $\sk[\nu]$ are close, the scheme presents stability areas nearly independent of $\utilde$. The case $\sk[e]=\sk[\nu]$ corresponds to the BGK scheme: the velocity field $\utilde$ does not play any role since (\ref{eq:mueq}) is verified.\\

Finally, for the d'Humières scheme, the results are independent of the choice of the moments. Indeed, the third order moment $X(X^2+Y^2)$ is equivalent to $\lambda^2 X+XY^2$ on the velocity network \cite{Fev:2014:1}: its relaxation is then equivalent to the relaxation of $XY^2$, $X$ being a conserved component. The same reasoning holds for the symmetrical moment. The moment $(X^2+Y^2)^2$ is equal to $X^2Y^2+\lambda^2(X^2+Y^2)$ on the velocity network. Its relaxation is equivalent to relax $X^2Y^2$ because $X^2+Y^2$ and $X^2Y^2$ are both in the eigenspace related to $\sk[e]$ for the TRT$_1$.\\

All these trends are independent of the direction of linearization $\theta$: the results are similar for $\theta=\pi/8,\pi/4,\pi/3$.\\

We now do the same job for the product equilibrium (\ref{eq:eqgeier}) restricting the study to the choice $\alpha=0$. Note that the combination of these moments and this equilibrium corresponds to the $\ddqn$ cascaded scheme \cite{Geier:2006:1,Fev:2014:0}. The TRT$_1$  (\ref{eq:trt1}), tuned by $\sk[e]$ and $\sk[\nu]$, and the direction $\theta=0$ are chosen. The table \ref{table:Maeqgeierdh1} is about the d'Humières scheme, the table \ref{table:Maeqgeierrec1} corresponds to the scheme relative to $\utilde=\vectV$.\\

\begin{table}
\centering\small
\grandtraittop
\begin{tabular}{@{}p{3.25cm}p{1cm}p{1cm}p{1cm}p{1.cm}p{1cm}p{1cm}p{1cm}p{1cm}@{}}
$n\qquad m$& 0&1&2&3&4&5&6&7\\ \hline
    0 & 0.42 & 0.42&0.39& 0.32&0.24&0.16&0.11&0.07 \\ 
     1 & 0.42 &0.42&0.41&0.38&0.31&0.20&0.14&0.09\\ 
     2 & 0.42 &0.42&0.41 &0.40&0.38&0.30&0.20&0.14\\ 
     3 &0.26&0.41&0.40&0.39&0.37 &0.32&0.28&0.20  \\ 
     4 &0.16&0.28&0.38&0.36&0.33&0.29&0.24&0.21\\
     5 &0.10&0.18&0.29&0.33 &0.31&0.26&0.21&0.19 \\
     6 &0.07&0.12&0.19&0.30&0.29&0.24&0.20&0.18\\ 
     7 &0.05&0.08&0.13&0.20&0.28&0.23&0.19&0.17  \\
\end{tabular}
\grandtraitbottom
\vspace{0.2cm}
 \caption{Highest stable $\vectV=(\Vx,0)$ in $\lambda$ units for the d'Humières scheme ($\utilde=\vectz$), with $\alpha=0$, of equilibrium (\ref{eq:eqgeier}). $\sk[e]=2-2^{-m}$ and $\sk[\nu]=2-2^{-n}$.}
 \label{table:Maeqgeierdh1}
\end{table}

\begin{table}
\centering\small
\grandtraittop
\begin{tabular}{@{}p{3.25cm}p{1cm}p{1cm}p{1cm}p{1.cm}p{1cm}p{1cm}p{1cm}p{1cm}@{}}
$n\qquad m$& 0&1&2&3&4&5&6&7\\ \hline
    0 & 0.42 & 0.42&0.42& 0.42&0.35&0.30&0.26&0.23 \\ 
     1 & 0.42 &0.42&0.42&0.41&0.39&0.35&0.32&0.28\\ 
     2 &0.42&0.41&0.41&0.40&0.40&0.35&0.31&0.29\\ 
     3 &0.41&0.41&0.40&0.39 &0.36&0.30&0.27&0.23  \\ 
     4 &0.40&0.40&0.39&0.36&0.33&0.28&0.23&0.20\\ 
     5 &0.35&0.37&0.36&0.33 &0.31&0.26&0.21&0.18 \\ 
     6 &0.31&0.33&0.33&0.31&0.29&0.25&0.20&0.17\\ 
     7 &0.28&0.30&0.31&0.29&0.27&0.24&0.19&0.17  \\
\end{tabular}
\grandtraitbottom
\vspace{0.2cm}
 \caption{Highest stable $\vectV=(\Vx,0)$ in $\lambda$ units for the scheme relative to $\utilde=\vectV$ with $\alpha=0$ of equilibrium (\ref{eq:eqgeier}). $\sk[e]=2-2^{-m}$ and $\sk[\nu]=2-2^{-n}$.}
 \label{table:Maeqgeierrec1}
\end{table}

The results are analogous to the ones associated with the equilibrium (\ref{eq:eqqian}).
The scheme relative to $\utilde=\vectV$ has bigger linear stability areas than the d'Humières scheme ($\utilde=\vectz$) when $\alpha=0$. The gain is more important when the relaxation parameters are far from each other. The velocity field impact is lightened when $\sk[e]$ and $\sk[\nu]$ are close.\\

We finally assess the influence of the equilibrium on the linear stability. Whatever the choice of $\utilde$ and $\vects$, the equilibrium (\ref{eq:eqgeier}) provides bigger stability areas than the truncated equilibrium (\ref{eq:eqqian}). Particularly, the BGK scheme associated with (\ref{eq:eqgeier}) is more stable than the one corresponding to (\ref{eq:eqqian}).\\

As a conclusion, the most important fact of the study is the following: the scheme relative to $\utilde=\vectV$ for $\alpha=0$ is more stable than for $\alpha=1$. Instead choosing a scheme relative to $\utilde=\vectV$ with a ``inappropriate'' set of moments can deteriorate the stability.

\subsection{Influence of the choice of the moments on the stability}\label{sub:stabmom}

The previous section has studied the stability of the relative velocity and d'Humières schemes for two choices of $\alpha$. This parameter seems to be crucial for the relative velocity schemes. The purpose of this section is to see more precisely its influence on the stability. It studies the stability properties of the schemes relative to $\utilde=\vectz$ and $\utilde=\vectV$ for a bigger range of $\alpha$. We show numerically and justify that $\alpha=0$ constitutes the better choice of moments.\\

We are interested in the stability of the relative velocity schemes for both sets of moments (\ref{eq:polal}) and (\ref{eq:polal2}): the discussion carries on the choice of the parameter $\alpha\in\R$ characterizing these moments. Both equilibria leading to the same trends, we focus on the truncated one (\ref{eq:eqqian})  linearized around $\vectV=(\Vx,0)\in\R^2$. Two sets of relaxation parameters $\vects$ are used: the TRT$_1$ (\ref{eq:trt1}) and the TRT$_2$ (\ref{eq:trt2}) where $\sk[e]=2-2^{-m}$ and $\sk[\nu]=\sk[p]=2-2^{-n}$ with $(m,n)=(0,3),(3,0),(0,7),(7,0),(7,7)$. The quantity (\ref{eq:maxVrho}) is drawn as a function of $\alpha$ in $[-1,1]$. A negative value of (\ref{eq:maxVrho}) means that the scheme is unstable for all $\vectV$.\\ 

\begin{figure}
  \begin{minipage}{0.45\linewidth} 
\begin{center}
\includegraphics[trim = 25mm 50mm 25mm 60mm,width=1.\textwidth]{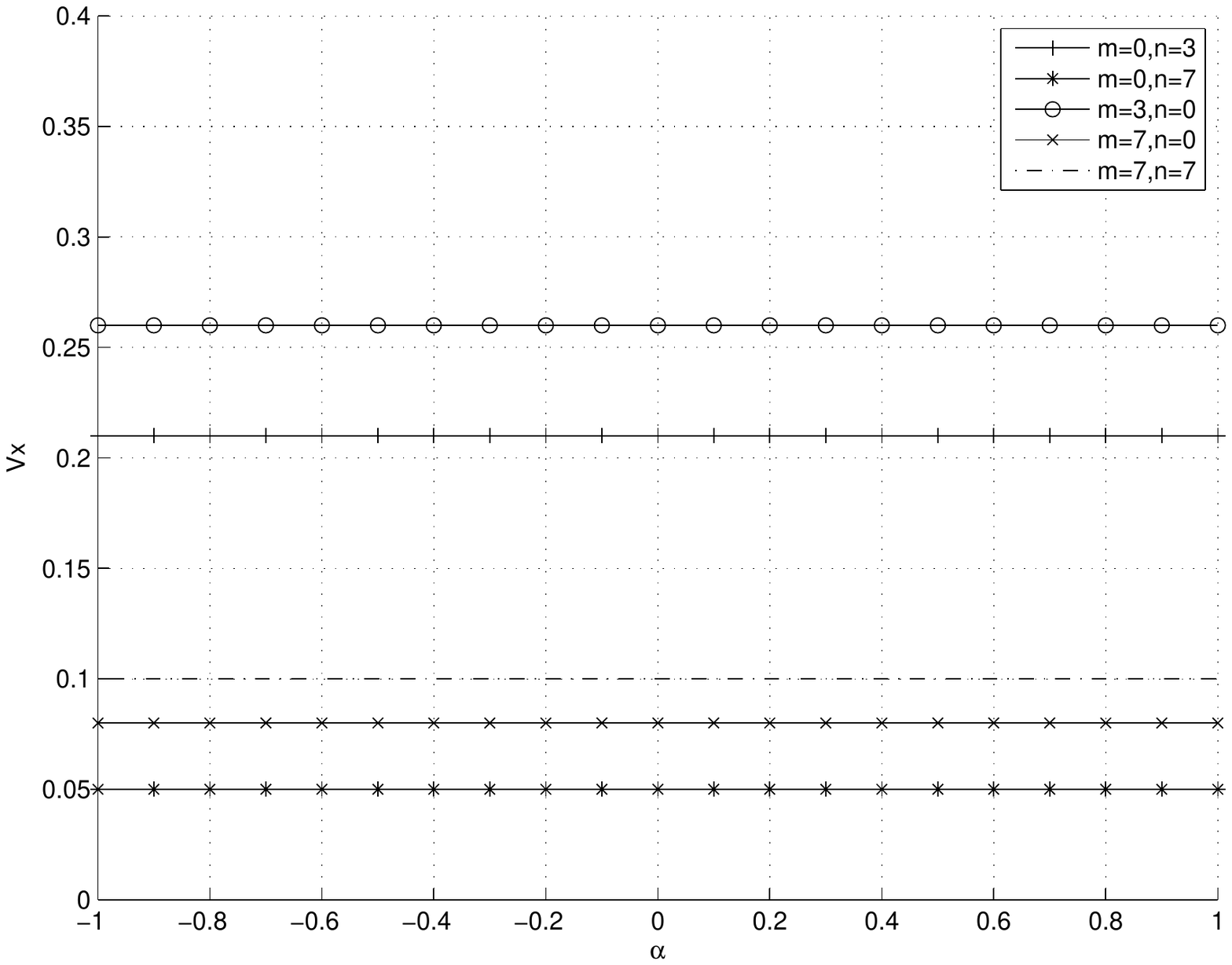}
\end{center}
\end{minipage}\hfill
\begin{minipage}{0.45\linewidth} 
\begin{center}
\includegraphics[trim = 25mm 50mm 25mm 60mm,width=1.\textwidth]{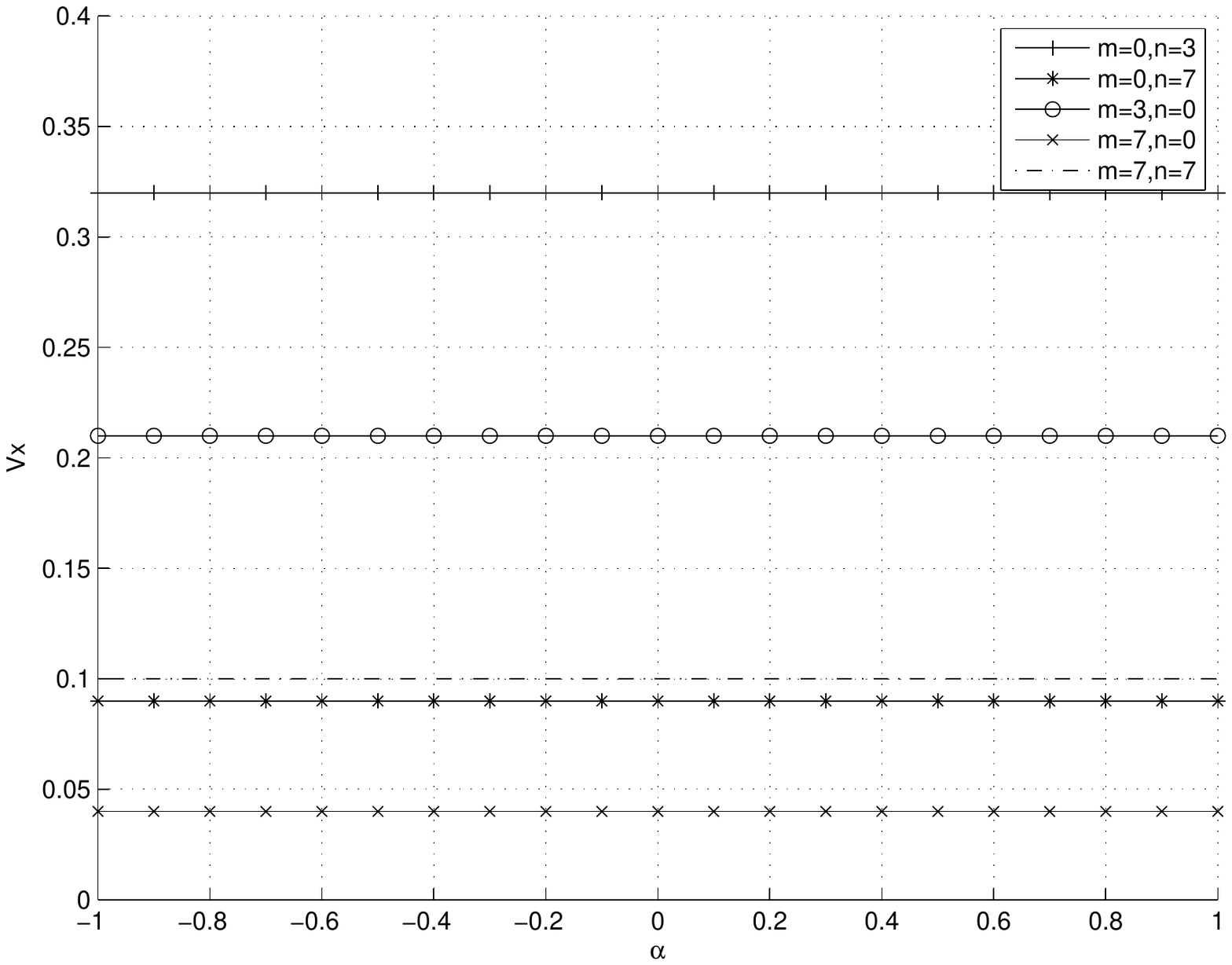}
\end{center}
\end{minipage}\\
\caption[$\ddqn$ $\utilde=\vectz$ : stabilité linéaire en $\Vx$ suivant le choix des moments 1]{Draw of $\Vx$ as a function of $\alpha$ for the d'Humières scheme with the moments  (\ref{eq:polal}). Left: TRT$_1$, $\sk[e]=2-2^{-m}$ and $\sk[\nu]=2-2^{-n}$.  Right: TRT$_2$, $\sk[e]=2-2^{-m}$ and $\sk[p]=2-2^{-n}$.}
\label{fig:stdhx3}
\end{figure}

\begin{figure}
  \begin{minipage}{0.45\linewidth} 
\begin{center}
\includegraphics[trim = 25mm 50mm 25mm 60mm,width=1.\textwidth]{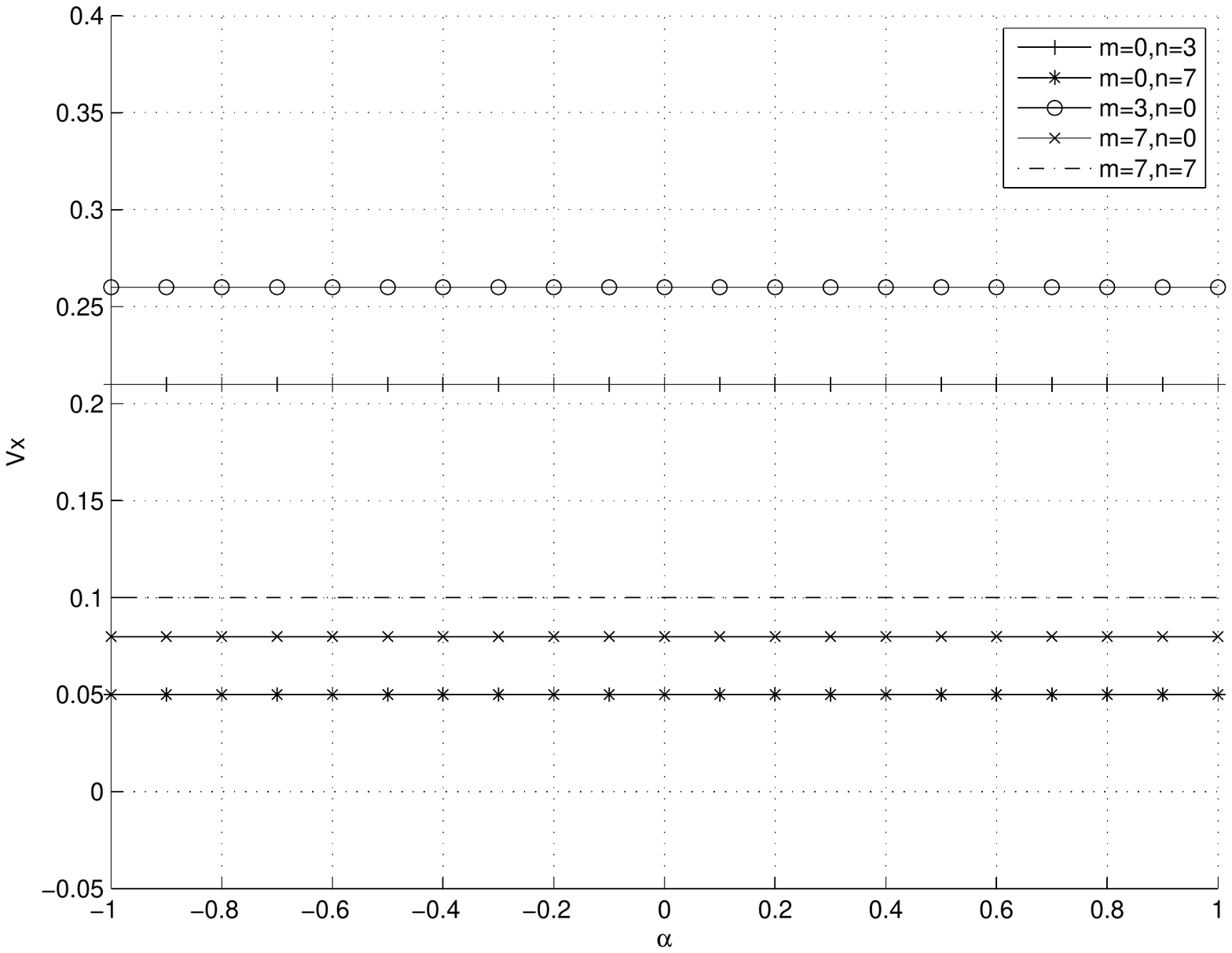}
\end{center}
\end{minipage}\hfill
\begin{minipage}{0.45\linewidth} 
\begin{center}
\includegraphics[trim = 25mm 50mm 25mm 60mm,width=1.\textwidth]{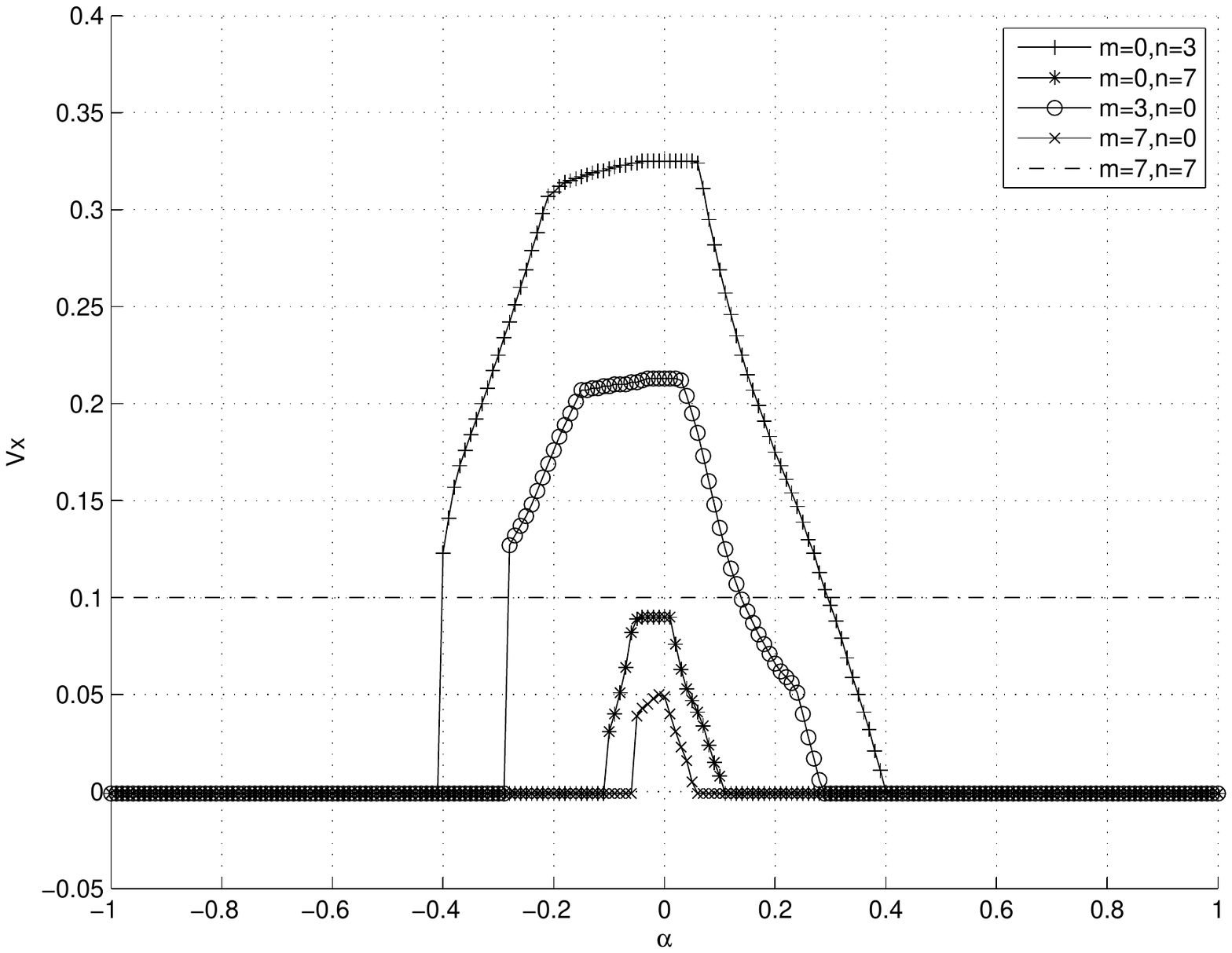}
\end{center}
\end{minipage}\\
\caption[$\ddqn$ $\utilde=\vectz$ : stabilité linéaire en $\Vx$ suivant le choix des moments 2]
{Draw of $\Vx$ as a function of $\alpha$ for the d'Humières scheme with the moments  (\ref{eq:polal2}). Left: TRT$_1$, $\sk[e]=2-2^{-m}$ and $\sk[\nu]=2-2^{-n}$.  Right: TRT$_2$, $\sk[e]=2-2^{-m}$ and $\sk[p]=2-2^{-n}$.}
\label{fig:stdhE}
\end{figure}

We first focus on the d'Humières scheme corresponding to $\utilde=\vectz$. The figures \ref{fig:stdhx3} and \ref{fig:stdhE} represent respectively the draws associated with the moments (\ref{eq:polal}) and (\ref{eq:polal2}). On each figure, the left draw is associated with the TRT$_1$ and the right one with the TRT$_2$.\\ 

For the moments (\ref{eq:polal}), the draws are independent of $\alpha$ whatever the TRT chosen and $\vects$. For the moments (\ref{eq:polal2}), the draw corresponding to the TRT$_1$ is independent of $\alpha$ unlike the TRT$_2$. The figure associated with the TRT$_2$ induces to choose $\alpha=0$: it corresponds to the maximum of the curve and the stability area decreases as $|\alpha|$ increases. As expected, the draw for $m=n=7$ corresponding to a BGK scheme is constant in $\alpha$. We notice that $\alpha=0$ belongs to the set of $\alpha$ maximizing the stability whatever the draw.\\

We can exhibit the origin of the dependence or independence on $\alpha$. Let's consider the moments (\ref{eq:polal}). For the d'Humières scheme, the relaxation of these moments is independent of $\alpha$. The last three moments of (\ref{eq:polal}) are $$\alpha X^3+XY^2, \alpha Y^3+X^2Y, \frac{\alpha}{2}(X^4+Y^4)+X^2Y^2.$$ 
Knowing that $X^3=\lambda^2 X$ on the velocity set \cite{Fev:2014:1}, the scheme is unchanged if we replace them by $$\lambda^2\alpha X+XY^2, \lambda^2\alpha Y+X^2Y, \frac{\lambda^2\alpha}{2}(X^2+Y^2)+X^2Y^2.$$
 Relaxing the moments (\ref{eq:polal}) is then equivalent to relax the same moments for $\alpha=0$. Indeed, $X$ and $Y$ are associated with some conserved moments and $X^2+Y^2$ has the same relaxation parameter $\sk[e]$ as the fourth order moment. It is thus consistent for this draws to be independent of $\alpha$.\\

We now focus on the moments (\ref{eq:polal2}): the parameter $\alpha$ appears only in the third order moments. For the TRT$_1$, $X^2+Y^2$ and the third order polynomials are relaxed with the same relaxation parameter $\sk[e]$. Choosing $$XY^2+\alpha(X^2+Y^2), X^2Y+\alpha(X^2+Y^2),$$ is then equivalent to choose $$XY^2, X^2Y,$$ and the scheme does not depend on $\alpha$ as the left draw of the figure \ref{fig:stdhE} shows it. For the TRT$_2$, $X^2+Y^2$ and the third order moments are relaxed with different relaxation parameters: it is expected to have a dependence on $\alpha$, excepted for the BGK case ($m=n=7$) involving only one relaxation parameter.\\
 
\begin{figure}
  \begin{minipage}{0.45\linewidth} 
\begin{center}
\includegraphics[trim = 25mm 50mm 25mm 60mm,width=1.\textwidth]{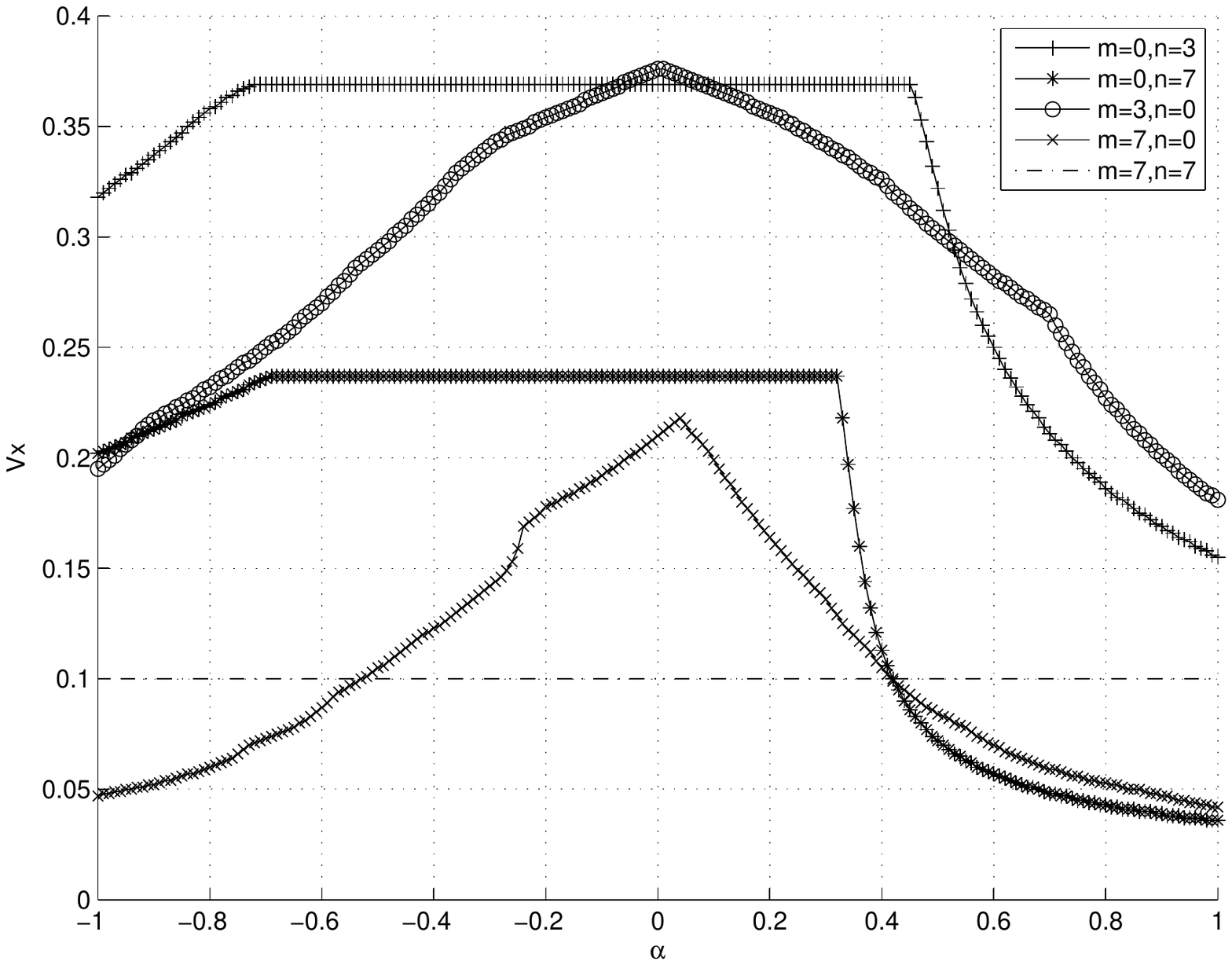}
\end{center}
\end{minipage}\hfill
\begin{minipage}{0.45\linewidth} 
\begin{center}
\includegraphics[trim = 25mm 50mm 25mm 60mm,width=1.\textwidth]{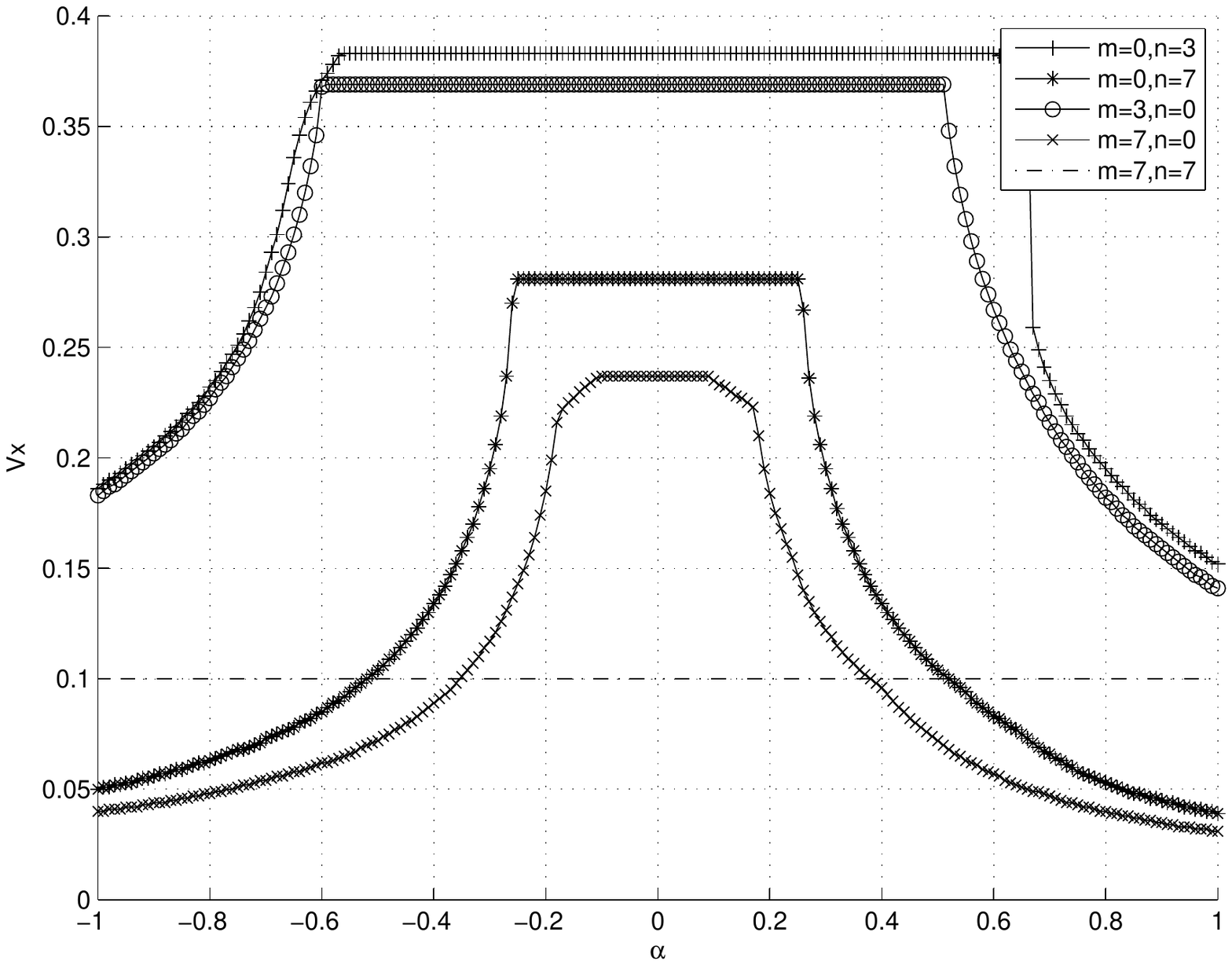}
\end{center}
\end{minipage}\\
\caption[$\ddqn$ $\utilde=\vectu$ : stabilité linéaire en $\Vx$ suivant le choix des moments 2]
{Draw of $\Vx$ as a function of $\alpha$ for the scheme relative to $\utilde=\vectV$ with the moments  (\ref{eq:polal}). Left: TRT$_1$, $\sk[e]=2-2^{-m}$ and $\sk[\nu]=2-2^{-n}$.  Right: TRT$_2$, $\sk[e]=2-2^{-m}$ and $\sk[p]=2-2^{-n}$.}
\label{fig:strecx3}
\end{figure}

\begin{figure}
  \begin{minipage}{0.45\linewidth} 
\begin{center}
\includegraphics[trim = 25mm 50mm 25mm 60mm,width=1.\textwidth]{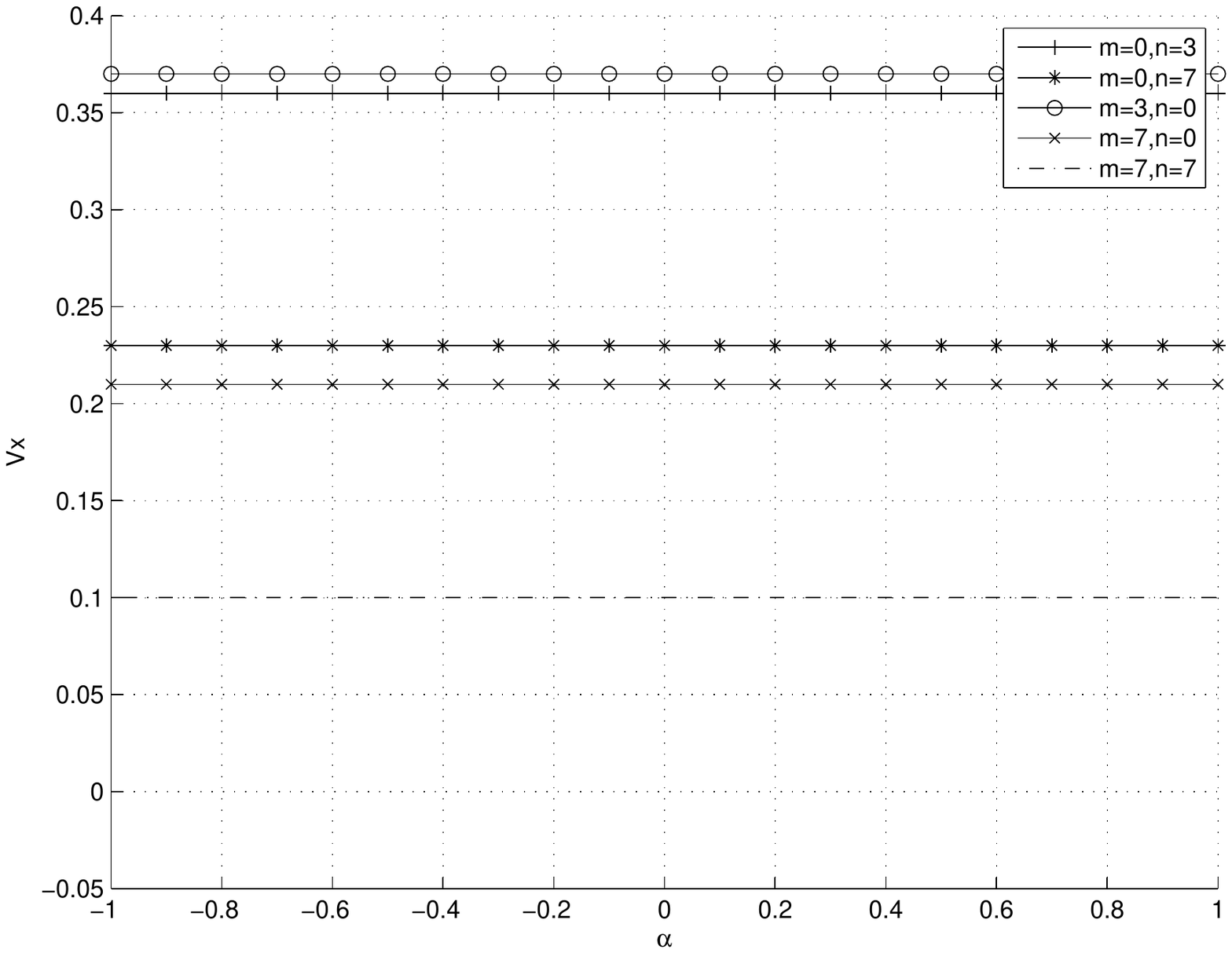}
\end{center}
\end{minipage}\hfill
\begin{minipage}{0.45\linewidth} 
\begin{center}
\includegraphics[trim =25mm 50mm 25mm 60mm,width=1.\textwidth]{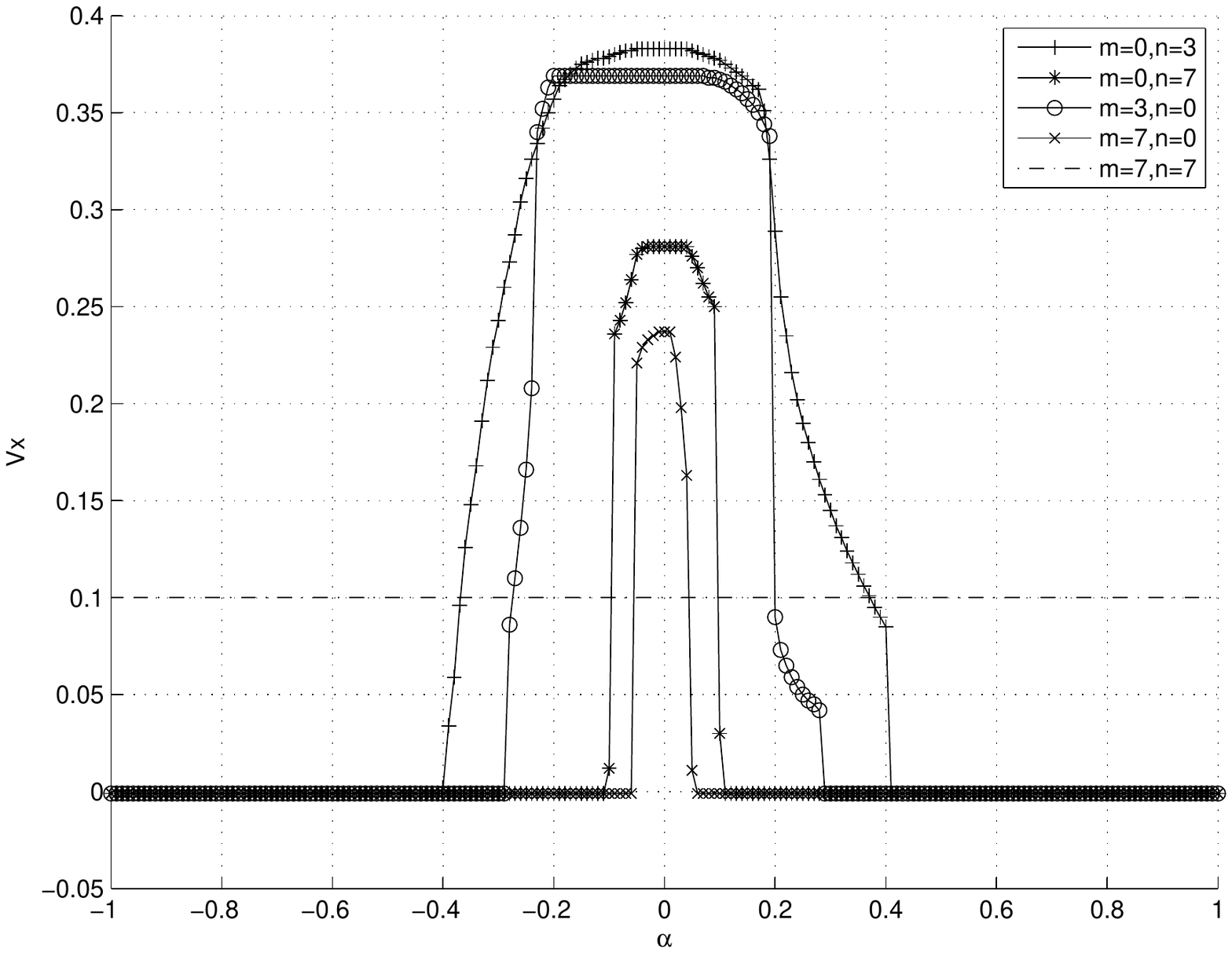}
\end{center}
\end{minipage}\\
\caption[$\ddqn$ $\utilde=\vectu$ : stabilité linéaire en $\Vx$ suivant le choix des moments 2]
{Draw of $\Vx$ as a function of $\alpha$ for the scheme relative to $\utilde=\vectV$ with the moments  (\ref{eq:polal2}). Left: TRT$_1$, $\sk[e]=2-2^{-m}$ and $\sk[\nu]=2-2^{-n}$.  Right: TRT$_2$, $\sk[e]=2-2^{-m}$ and $\sk[p]=2-2^{-n}$.}
\label{fig:strecE}
\end{figure}

 We now do the same job for the scheme relative to $\utilde=\vectV$. The figure \ref{fig:strecx3} is associated with the moments (\ref{eq:polal}) and the figure \ref{fig:strecE} with the moments (\ref{eq:polal2}).\\

The stability of the scheme relative to $\utilde=\vectV$ depends on $\alpha$ whatever the moments (figure \ref{fig:strecx3}). The maximum is reached for $\alpha=0$ whatever the choice of $\vects$. For the moments (\ref{eq:polal2}), the stability of the TRT$_1$ is not linked to $\alpha$ (figure \ref{fig:strecE} on the left side). 
Instead, this parameter is influential for the TRT$_2$ (figure \ref{fig:strecE} on the right side): $\alpha=0$ still corresponds to the optimum.\\
 
 We interpret the figure \ref{fig:strecx3} corresponding to the moments (\ref{eq:polal}). Because $X^3=\lambda^2 X$ and $Y^3=\lambda^2 Y$ on the velocity set, relaxing the relative moments associated with (\ref{eq:polal}) is equivalent to relax 
\begin{equation*}\label{eq:momeqalpha_4}
1,X,Y,X^2+Y^2,X^2-Y^2,XY,\overline{\Pk[6]}(\utilde,\alpha),\overline{\Pk[7]}(\utilde,\alpha),\overline{\Pk[8]}(\utilde,\alpha),
\end{equation*}
where
\begin{align}
\overline{\Pk[6]}(\utilde,\alpha)&=XY^2+\alpha(-3\utx X^2+(\lambda^2-3(\utx)^2)X+\utx(\lambda^2-(\utx)^2)),\label{eq:P6}\\
\overline{\Pk[7]}(\utilde,\alpha)&=X^2Y+\alpha(-3\uty Y^2+(\lambda^2-3(\uty)^2)Y+\uty(\lambda^2-(\uty)^2)),\nonumber
\end{align}
and
\begin{multline*}
\overline{\Pk[8]}(\utilde,\alpha)=X^2Y^2+\frac{\alpha}{2}\Big((\lambda^2+6(\utx)^2)X^2+(\lambda^2+6(\uty)^2)Y^2\\
+2\utx(-\lambda^2+4(\utx)^2)X+2\uty(-\lambda^2+4(\uty)^2)Y\\
-3(\utx)^2(\lambda^2-(\utx)^2)-3(\uty)^2(\lambda^2-(\uty)^2)\Big).
\end{multline*}
Let's observe the equivalent class of the third order moment given by (\ref{eq:P6}).
The dependence on $\alpha$ of the stability comes from the term $-3\alpha\utx X^2$. Indeed, relaxing (\ref{eq:P6}) is equivalent to relax $XY^2-3\alpha\utx X^2$ since the moments corresponding to the polynomials $1$ and $X$ are conserved by the collision. On the contrary, the moment associated with $X^2$ is not conserved. For the TRT$_1$, it is a linear combination of the moments $X^2+Y^2$ and $X^2-Y^2$ associated with different relaxation parameters $\sk[e]$ and $\sk[\nu]$. For the TRT$_2$, it is associated with $\sk[e]$ whereas $\overline{\Pk[6]}$ corresponds to $\sk[p]$.\\

These remarks justify the introduction of the moments (\ref{eq:polal2}) to study the influence of the non conserved components $X^2$ and $Y^2$. The figure \ref{fig:strecE} implying the moments (\ref{eq:polal2}) gives similar results as its analogous for $\utilde=\vectz$ (figure \ref{fig:stdhE}): the same interpretation is still valid. Note that for $\alpha=0$, the areas are bigger with $\utilde=\vectV$ (figure \ref{fig:strecE}) than with $\utilde=\vectz$ (figure \ref{fig:stdhE}). This confirms the observations of the section \ref{sub:stlindH}.

\section{Stability for the Kelvin-Helmholtz instability}\label{se:KHstab}

The purpose of this section is to confirm on a non linear test case the previous linear stability results: this test case is the Kelvin-Helmhotz instability \cite{Minion:1997:0,Dellar:2001:0}.\\

We compare six versions of the relative velocity $\ddqn$ scheme to study the influence of the moments, of the velocity field $\utilde$ and of the equilibrium. We consider the scheme associated with $\alpha=0$ relative to $\utilde=\vectz$ and $\utilde=\vectu$ (the fluid velocity) for the equilibria (\ref{eq:eqqian}) and (\ref{eq:eqgeier}). We compare it to the choice $\alpha=1$ for the relative velocities $\vectz$ and $\vectu$ with the equilibrium (\ref{eq:eqqian}). We choose not to consider the product equilibrium (\ref{eq:eqgeier}) for $\alpha=1$, this equilibrium being introduced for the moments of the cascaded scheme \cite{Geier:2006:1}. We work with the TRT$_1$ defined by (\ref{eq:trt1}): unless otherwise specified, $\sk[e]$ et $\sk[\nu]$ are fixed by
\begin{equation*}
 \mu=\frac{\lambda^2\Delta\sig[e]}{3},\quad\nu=\frac{\lambda^2\Delta\sig[\nu]}{3},
 \end{equation*}
 where $\sig[e]=1/\sk[e]-1/2$ and $\sig[\nu]=1/\sk[\nu]-1/2$, so that the viscosities $\mu$ and $\nu$ are set to $0.0366$ and $10^{-4}$.\\

We test the stability of the scheme by increasing the velocity $U$ defining the initial shear layers  \begin{equation*}\label{eq:KHCI_4}
 u^x(x,y,0)=\left\{\begin{split}
 U\tanh(k(y-\tfrac{1}{4}))&~ {\rm if}~ y\leq\tfrac{1}{2}\\
 U\tanh(k(\tfrac{3}{4}-y))&~ {\rm if}~ y>\tfrac{1}{2}
 \end{split}
 \right., \quad (x,y)\in[0,1]^2,
 \end{equation*} 
  \begin{equation*}\label{eq:KHCI2_4}
 u^y(x,y,0)=U\delta\sin(2\pi(x+\tfrac{1}{4})),\quad (x,y)\in[0,1]^2.
 \end{equation*} 
 This velocity $U$ is chosen as Ma$/\sqrt{3}$ for Ma$\in\R$ the Mach number.  The parameters $k$ and $\delta$ controlling the width of the shear layers and the magnitude of the initial data are set to $80$ and $0.05$.\\

\begin{figure}[!t]\label{fig:nu}
\begin{center}
\includegraphics[trim = 0mm 30mm 0mm 30mm,width=1.\textwidth]{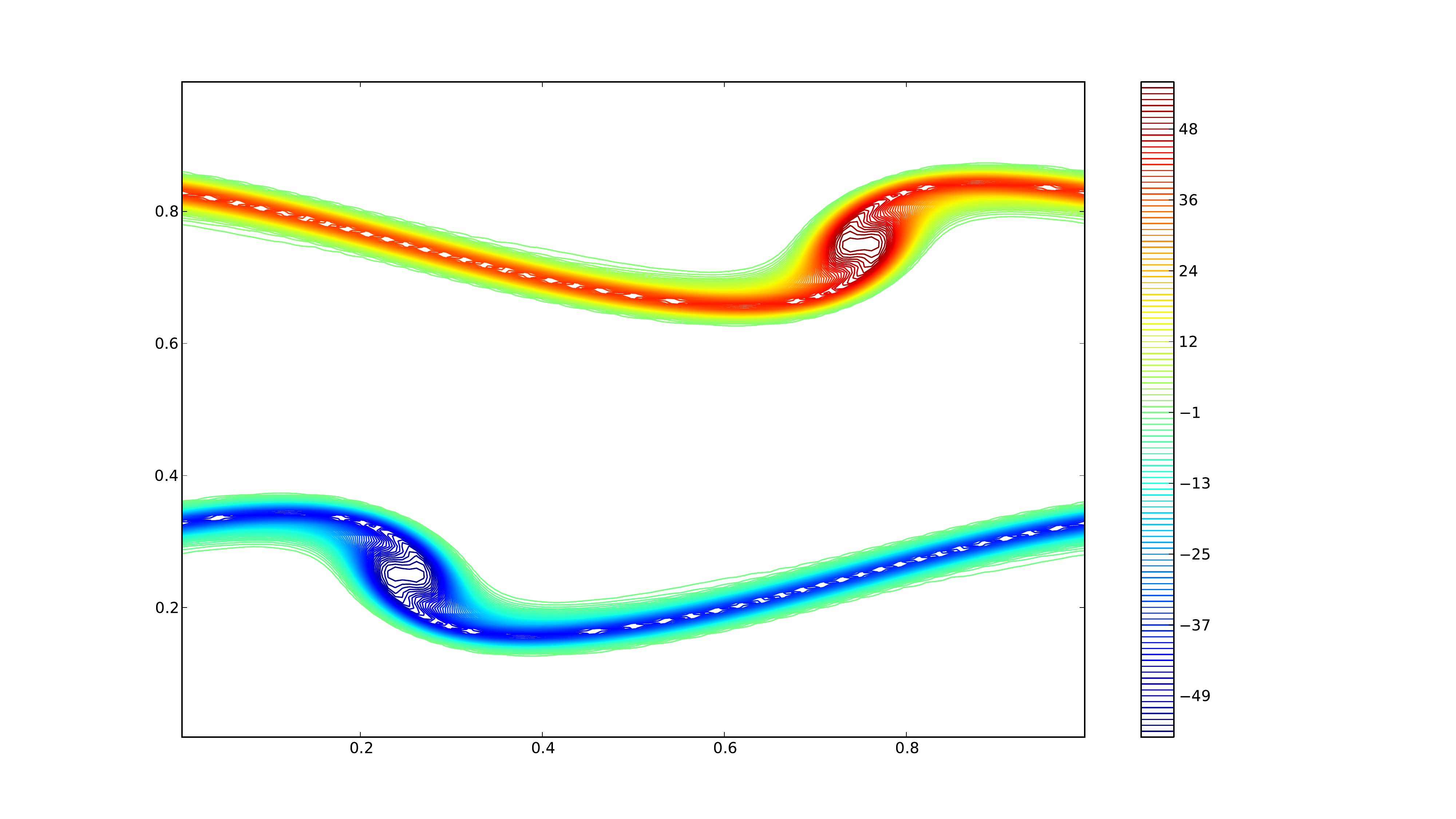}
\end{center}
\vspace{0.2cm}
\caption{Vorticity draw at t=0.6.}
\label{fig:vort1}
\end{figure}

\begin{figure}[!t]\label{fig:nu}
\begin{center}
\includegraphics[trim = 0mm 30mm 0mm 30mm,width=1.\textwidth]{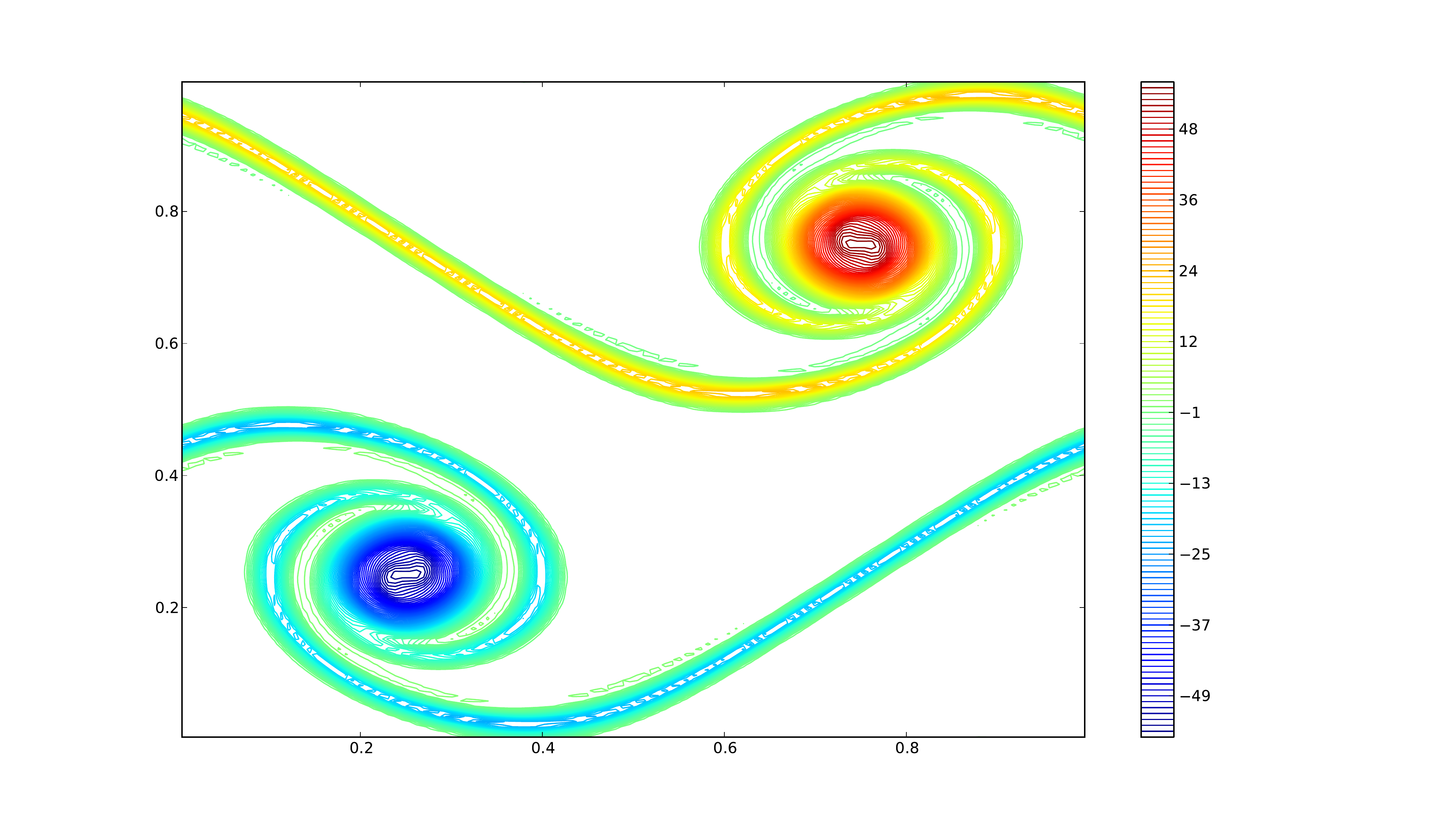}
\end{center}
\vspace{0.2cm}
\caption{Vorticity draw at à t=1.}
\label{fig:vort2}
\end{figure}

We first validate the vorticity draws obtained in \cite{Ricot:2009:0,Dellar:2001:0,Minion:1997:0} using the scheme relative to the fluid velocity $\vectu$ for the second order truncated equilibrium (\ref{eq:eqqian}). This vorticity is defined by
 $$\omega=\Dx u_y-\partial_y u_x.$$
 For this simulation, the domain is constituted of $128\times128$ points, the Mach number is fixed at $0.04$ ($\lambda$ is chosen as in \cite{Dellar:2001:0} so that $U=1$). The figures \ref{fig:vort1} and \ref{fig:vort2} are the vorticity plots at time $t=0.6$ and $t=1$.\\


We now present a stability analysis depending on the different parameters for $\lambda=1$.
We expect to confirm the linear stability results. The scheme is considered stable if it has not broken after 2000 iterations. The table \ref{table:Madx} contains the maximal stable Mach number Ma for different meshes at $0.01$ close. The table $\ref{table:Re}$ presents the greater Reynolds number $Re=1/\nu$ stable at $1000$ close for different meshes and Ma$=0.09$. Since we discuss on the Reynolds number, the viscosity $\nu$ becomes a parameter.\\


\begin{table}
\centering\small
\grandtraittop
\begin{tabular}{@{}p{4.65cm}p{1cm}p{1cm}p{1cm}p{1cm}p{1cm}p{1cm}p{1cm}@{}}
 Space step $\dx$ &$1/16$&$1/32$&$1/64$&$1/128$&$1/256$&$1/512$&$1/1024$\\
  Corresponding $\sk[e]$  & 0.44 & 0.25  &0.13 &0.07&0.03& 0.02&0.01\\ 
  Corresponding $\sk[\nu]$ & 1.98 & 1.96  &1.93 &1.86&1.73& 1.53&1.24\\ \hline
   $\alpha=0$, $\utilde=\vectz$, equilibrium (\ref{eq:eqqian}) &0.18&0.13&0.12 &0.12&0.12&0.12&0.12\\ 
     $\alpha=0$, $\utilde=\vectu$, equilibrium (\ref{eq:eqqian})&0.96&0.82&0.62&0.49&0.43&0.39&0.39\\ 
     $\alpha=0$, $\utilde=\vectz$, equilibrium (\ref{eq:eqgeier})&0.18&0.13&0.12&0.12&0.12&0.12&0.12  \\ 
     $\alpha=0$, $\utilde=\vectu$, equilibrium (\ref{eq:eqgeier})  &0.92&0.80&0.62&0.50&0.43&0.39&0.39 \\ 
     $\alpha=1$, $\utilde=\vectz$, equilibrium (\ref{eq:eqqian}) &0.18&0.13&0.12 &0.12&0.12&0.12&0.12 \\ 
     $\alpha=1$, $\utilde=\vectu$, equilibrium (\ref{eq:eqqian})&0.09&0.07&0.06& 0.05&0.05&0.05&0.05  \\ 
\end{tabular}
\grandtraitbottom
\vspace{0.2cm}
 \caption{Maximum of Ma stable according to the mesh. The two last columns indicates a convergence of each scheme as $\dx$ decreases.}
  \label{table:Madx}
\end{table}

%

\begin{table}
\centering\small
\grandtraittop
\begin{tabular}{@{}p{6.7cm}p{1.5cm}p{1.5cm}p{1.5cm}p{1.5cm}@{}}
Space step $\dx$ &$1/16$&$1/32$&$1/64$&$1/128$\\
  Corresponding $\sk[e]$ & 0.44 & 0.25  &0.13 &0.07  \\ 
  Corresponding $\sk[\nu]$ & 1.98 & 1.96  &1.93 &1.86 \\ \hline
     $\alpha=0$, $\utilde=\vectz$, equilibrium (\ref{eq:eqqian}) &$\sk[\nu]=2$&$21.10^3$&$17.10^3$&$17.10^3$\\ 
     $\alpha=0$, $\utilde=\vectu$, equilibrium (\ref{eq:eqqian})&$\sk[\nu]=2$&$\sk[\nu]=2$&$\sk[\nu]=2$&$\sk[\nu]=2$\\ 
     $\alpha=0$, $\utilde=\vectz$, equilibrium (\ref{eq:eqgeier})&$\sk[\nu]=2$&$21.10^3$&$17.10^3$&$17.10^3$\\ 
     $\alpha=0$, $\utilde=\vectu$, equilibrium (\ref{eq:eqgeier}) &$\sk[\nu]=2$&$\sk[\nu]=2$&$\sk[\nu]=2$&$\sk[\nu]=2$\\ 
    $\alpha=1$, $\utilde=\vectz$, equilibrium (\ref{eq:eqqian}) &$\sk[\nu]=2$&$21.10^3$&$17.10^3$&$17.10^3$\\ 
    $\alpha=1$, $\utilde=\vectu$, equilibrium (\ref{eq:eqqian})&$10.10^3$&$6.10^3$&$4.10^3$&$4.10^3$\\ 
\end{tabular}
\grandtraitbottom
\vspace{0.2cm}
 \caption{Maximum of the Reynolds number stable for ${\rm Ma}=0.09$ according to the mesh.}
  \label{table:Re}
\end{table}

We obtain results consistent with the linear stability study. First, choosing a scheme relative to $\utilde=\vectu$ has a positive effect if $\alpha=0$, negative if $\alpha=1$. We must choose the moments of the $\ddqn$ cascaded scheme to improve the stability. This improvement occurs whatever the equilibrium and the mesh: the stability limit $\sk[\nu]=2$ is stable (table \ref{table:Re}) and high Mach numbers are reached for this scheme (table \ref{table:Madx}).
Second, the d'Humières scheme is independent of $\alpha$ as for the linear stability study. Its stability area is smaller than the scheme relative to $\utilde=\vectu$ when $\alpha=0$, greater when $\alpha=1$. Finally, the equilibrium does not influence a lot the stability unlike the linear case. The obtained values are close whatever the choice of the equilibrium. It is important to note that the table \ref{table:Madx} exhibits a convergence of all the schemes as $\dx$ decreases.\\

We now characterize the behaviour of the scheme when the diffusion is weak (when the relaxation parameters are close to $2$). The table \ref{table:Mamu} presents the maximal Ma stable for decreasing bulk viscosity $\mu$. The domain is constituted of $128^2$ points and $\nu$ is still equal to $10^{-4}$.\\


\begin{table}
\centering\small
\grandtraittop
\begin{tabular}{@{}p{5.8cm}p{1.3cm}p{1.3cm}p{1.3cm}p{1.3cm}p{1.3cm}@{}}
  Viscosity $\mu$ & $10^{-2}$ &   $0.5.10^{-2}$  & $10^{-3}$  & $10^{-4}$  & $10^{-5}$ \\ 
     Corresponding $\sk[e]$  &0.23&0.41&1.13&1.86&1.98 \\ 
     Corresponding $\sk[\nu]$ &1.86&1.86&1.86&1.86&1.86 \\ \hline
    $\alpha=0$, $\utilde=\vectz$, equilibrium (\ref{eq:eqqian}) &0.22&0.29&0.43&0.38&0.32\\ 
     $\alpha=0$, $\utilde=\vectu$, equilibrium (\ref{eq:eqqian})&0.73&0.76&0.68&0.63&0.60\\ 
     $\alpha=0$, $\utilde=\vectz$, equilibrium (\ref{eq:eqgeier})& 0.22&0.30&0.45&0.38&0.32\\ 
     $\alpha=0$, $\utilde=\vectu$, equilibrium (\ref{eq:eqgeier})  &0.72&0.76&0.76&0.63&0.61\\ 
     $\alpha=1$, $\utilde=\vectz$, equilibrium (\ref{eq:eqqian}) &0.22&0.29&0.43&0.38&0.32\\ 
     $\alpha=1$, $\utilde=\vectu$, equilibrium (\ref{eq:eqqian})&0.10&0.14&0.32&0.38&0.32 \\ 

\end{tabular}
\grandtraitbottom
\vspace{0.2cm}
 \caption{Maximum of Ma stable according to $\mu$.}
  \label{table:Mamu}
\end{table}


This table is also consistent with the linear stability study. When $\sk[e]$ and $\sk[\nu]$ are far from each other, the linear case (the tables \ref{table:Mageier} and \ref{table:Maeqgeierrec1}) presents an important gain for the scheme relative to $\utilde=\vectu$ for $\alpha=0$ whatever the equilibrium. These results are confirmed by the first columns of the table \ref{table:Mamu} corresponding to take a big bulk viscosity. When $\mu$ tends to $0$, the stability areas for the different choices of $\utilde$ are expected to be close at fixed equilibrium: this case corresponds to close parameters $\sk[e]$ and $\sk[\nu]$, regime where the linear stability results are homogeneous in $\utilde$. This behaviour is confirmed by the table \ref{table:Mamu}: indeed, the four cases associated with the equilibrium (\ref{eq:eqqian}) have the same stability areas when the bulk viscosity is smaller than $10^{-4}$. Similarly, the two cases associated with the equilibrium (\ref{eq:eqgeier}) have close stability areas for these viscosities.\\

The table \ref{table:Marec} deals with the influence of the velocity field $\utilde$ on the stability: other choices than $\vectz$ and $\vectu$ are considered. We determine the maximal Mach number stable for different $\utilde$ according to the choice of the moments. We study the two choices $\alpha=0$ and $\alpha=1$ for a mesh of $128^2$ points.\\


\begin{table}
\centering\small
\grandtraittop
\begin{tabular}{@{}p{3.7cm}p{0.95cm}p{0.95cm}p{0.95cm}p{0.95cm}p{0.95cm}p{0.95cm}p{0.95cm}p{0.95cm}@{}}
   $\utilde$ & $\vectz$& $0.2\vectu$& $0.4\vectu$ & $0.6\vectu$ & $0.8\vectu$ & $\vectu$ & $1.2\vectu$&$1.4\vectu$ \\ \hline
     $\alpha=0$, equilibrium (\ref{eq:eqqian})   & 0.12 &0.15&0.21&0.34&0.60&0.49&0.42&0.33\\ 
     $\alpha=1$,  equilibrium (\ref{eq:eqqian})  &0.12&0.11&0.09&0.07& 0.06&0.05& 0.05 &0.04\\ 
    $\alpha=0$, equilibrium (\ref{eq:eqgeier})   & 0.12 &0.15&0.21&0.34&0.60&0.50&0.42&0.33\\

\end{tabular}
\grandtraitbottom
\vspace{0.2cm}
 \caption{Maximum of Ma stable according to $\utilde$.}
  \label{table:Marec}
\end{table}


This table is an evidence of the importance of the moments for the relative velocity schemes. Taking a velocity different from $\vectz$ provides stability improvements only for $\alpha=0$. These moments improve the numerical stability for $\utilde=\vectu$ compared to the d'Humières scheme whatever the equilibrium. Instead, choosing $\utilde\neq\vectz$ for the moments (\ref{eq:polddqn}) deteriorates the stability of the scheme. The most stable choice for $\alpha=1$ corresponds to the d'Humières scheme.

\section{Conclusion}

We have studied the numerical stability of the relative velocity $\ddqn$ scheme with two conservation laws. A linear stability study was presented and strenghtened by a non linear test case for the compressible Navier-Stokes equations: the Kelvin-Helmholtz instability. The main conclusion of the article is the following: the relative velocity schemes improve or deteriorate the stability of the d'Humières schemes and it depends strongly on the choice of the moments.  An improvement occurs if the moments of the cascaded scheme are chosen whatever the equilibrium. It is bigger when one viscosity is very small and the other is important. The usual set of moments and its orthogonalized analogous deteriorates the stability of the d'Humières scheme. This degradation originates from the presence of second order components in the third and fourth order moments. These components do not appear for the moments of the cascaded scheme that explains the better stability behaviour.

\bibliographystyle{plain}
\bibliography{Bibliographie}

\end{document}